%% file: diamondfree-FINAL-11-2012.tex
\documentclass[preprint, letter, final]{elsarticle}
\input{diamond-macro}

\begin{document}
\title{On diamond-free subposets of the Boolean lattice}
\author{Lucas~Kramer}\ead{ljkramer@iastate.edu}

\author{Ryan~R.~Martin\fnref{fn1}}\ead{rymartin@iastate.edu}

\author{Michael Young\fnref{fn2}}\ead{myoung@iastate.edu}

\address{Department of Mathematics, Iowa State University, Ames, Iowa 50011}

\fntext[fn1]{Corresponding author. This author's research partially supported by NSF grant DMS-0901008.}
\fntext[fn2]{This author's research partially supported by a fellowship funded by NSF grant DMS-0946431.}
\begin{keyword}
   forbidden subposets \sep extremal set theory \sep diamond-free \sep flag algebras
\end{keyword}

\begin{abstract}
The Boolean lattice of dimension two, also known as the diamond, consists of four distinct elements with the following property: $A\subset B,C\subset D$.  A diamond-free family in the $n$-dimensional Boolean lattice is a subposet such that no four elements form a diamond. Note that elements $B$ and $C$ may or may not be related.

There is a diamond-free family in the $n$-dimensional Boolean lattice of size $(2-o(1)){n\choose\lfloor n/2\rfloor}$. In this paper, we prove that any diamond-free family in the $n$-dimensional Boolean lattice has size at most $(2.25+o(1)){n\choose\lfloor n/2\rfloor}$.  Furthermore, we show that the so-called Lubell function of a diamond-free family in the $n$-dimensional Boolean lattice which contains the empty set is at most $2.25+o(1)$, which is asymptotically best possible.

AMS 2010 Subject Classification: Primary 06A07; Secondary 05D05, 05C35
\end{abstract}
\maketitle
\section{Introduction and Background}

In 1928, Sperner~\cite{Sp:1928} proved that the largest family of subsets of $[n]\stackrel{\rmdef}{=}\{1,2,\ldots,n\}$ with the property that no one of these sets contains another has size $\nntwo$. Sperner's Theorem can be derived from the YBLM (or LYM) inequality which states that if $\mathcal{A}$ is an antichain in $\booln$, then $\sum_{A\in\mathcal{A}}\binom{n}{|A|}^{-1}\leq 1$.  See, for example, De Bonis and Katona~\cite{BoKa:2007}.

The $n$-dimensional Boolean lattice, $\booln$, denotes the partially ordered set (or poset) $(2^{[n]},\subseteq)$, where, for every finite set $S$, $2^S$ denotes the set of subsets of $S$. (The poset $\booln$ is also called $\mathcal{Q}_n$ in other literature, such as~\cite{AxMaMa:2011}.)

We say that poset $\PP'=(\PP',\leq')$ is a \textit{(weak) subposet} of poset $\PP=(\PP,\leq)$ if there exists an injective function $f:P'\rightarrow P$ that preserves the partial ordering.  That is, if $u\leq' v$ in $\PP'$ then $f(u)\leq f(v)$ in $\PP$.  See Stanley~\cite{St:1997}.

Consider a family $\F$ of members of $2^{[n]}$. Then $\F$ can be viewed as a subposet of $\booln$ via the inherited relation. If $\F$ contains no  $\PP$ as a (weak) subposet we say $\F$ is \textit{$\PP$-free}. As in much of the literature (e.g.,~\cite{GrLiLu:2012}), we denote $\La (n,\PP)$ to be the size of the largest family $\F$ of elements of $\booln$ that is $\PP$-free.

A \textit{layer} of $\booln$ is a family of sets of the same cardinality, denoted $\binom{[n]}{k}$ for some integer $k$. A \textit{chain} in a poset is a set of elements, all of which are pairwise related. A chain with $k$ elements is said to have \textit{length $k-1$} and is denoted $P_k$.  So, Sperner's theorem says that a subposet in $\booln$ with no $\PP_2$ has size at most $\nntwo$.  Erd\H{o}s~\cite{Er:1945} generalized this result to show that, for any $k\geq 2$, a $\PP_k$-free family in $\booln$ has size at most $(k-1)\nntwo$.  Furthermore, the size is exactly the sum of the $k-1$ largest layers.

For many other subposets, $\PP$, the size of the largest $\PP$-free family in the $n$-dimensional Boolean lattice is not known, even asymptotically.  The simplest example, the diamond poset, consisting of only four elements, still has a wide gap between the largest known families and a proven upper bound.  In this paper, we prove that every diamond-free family in the Boolean lattice has size at most $(2.25+o(1))\nntwo$.  The lower bound of $(2-o(1))\nntwo$ is found by simply taking the layers $\binom{[n]}{\lfloor n/2\rfloor}$ and $\binom{[n]}{\lfloor n/2\rfloor+1}$.



 If $n$ is a fixed integer, let $\mathcal{B}(n,k)$ denote the collection of  subsets of the $k$ middle sizes, that is, the collection of sets of size $\left\lfloor (n-k+1)/2 \right\rfloor,\cdots,\left\lfloor (n+k-1)/2 \right\rfloor$ or $\left\lceil (n-k+1)/2 \right\rceil,\cdots,\left\lceil (n+k-1)/2 \right\rceil$ depending on the parity of $n$. If $n$ is a fixed integer and $\sum (n,k)$ denotes the sum of the $k$ largest binomial coefficients of the form ${n\choose l}$, then if $\PP_k$ denotes a chain with $k$ elements, Erd\H{o}s' result is as follows:

\begin{thm}[Erd\H{o}s~\cite{Er:1945}]\label{thm:Erdos}
For $n\geq k-1\geq 1$, $\La (n,\PP_k)=\sum (n,k-1)$. Moreover, the $\PP_k$-free families of maximum size in $\booln$ are given by $\mathcal{B}(n,k)$.
\end{thm}



Asymptotically, Theorem~\ref{thm:Erdos} gives that for any fixed $k$, $\La (n,\PP_k)=(k-1-o(1))\nntwo$. There are a number of results on the problem of finding $\La(n,\PP)$ for several other posets $\PP$. In particular, there are results regarding asymptotic bounds on $\La(n,\PP)/\nntwo$.  A number of these results are catalogued in the work of Griggs, Li and Lu~\cite{GrLiLu:2012}.  We will highlight a few here:

Let $\mathcal{V}_r$ be the $r$-fork poset defined to be the poset that has elements $A<B_1,B_2,\ldots, B_r$ for $r\geq 2$. Katona and Tarj\'an \cite{KaTa:1983} began studying these for $r=2$ and obtained bounds on $\La (n,\mathcal{V}_2)$ and Katona and DeBonis \cite{BoKa:2007} expanded these bounds for general $V_r$ for $r\geq 2$. They showed that $$ \left(1+\frac{r-1}{n} +\Omega\left(\frac{1}{n^2}\right)\right)  \leq \frac{\La (n,\mathcal{V}_r)}{\nntwo} \leq \left(1+2\frac{r-1}{n} +O\left(\frac{1}{n^2}\right)\right) .$$
From this we can see that $\La(n,\mathcal{V}_r) \sim  \nntwo.$




There have been a number of results that establish the asymptotic value of $\La(n,\PP)/\nntwo$ for certain posets $\PP$.  See, for instance, \cite{Th:1998,GrLu:2009,BoKaSw:2005,BoKa:2007}.  In each case, the asymptotic value of $\La(n,\PP)/\nntwo$ is determined and was found to be an integer.  As a result, Griggs and Lu~\cite{GrLu:2009} stated the following conjecture:

\begin{conj}\label{conj:GLL}
For every finite poset $\PP$, the limit $\pi (\PP)\stackrel{\rmdef}{=}\lim\limits_{n\rightarrow \infty} \La(n,\PP)/\nntwo$ exists and is an integer.
\end{conj}

This conjecture was verified by Griggs and Lu~\cite{GrLu:2009} for tree posets of height 2 as well as for the \textit{crown poset} $\mathcal{O}_{2k}, k\geq 2$, when $k$ is even and $k\geq 4$.  The crown is the poset of height 2 that is a cycle of length $2k$ in the Hasse diagram. However, when $k$ is odd, Griggs and Lu were only able to show that $\mathcal{O}_{2k}/\nntwo$ is at most $1+\frac{1}{\sqrt{2}}+o(1)$.


Griggs, Li and Lu~\cite{GrLiLu:2012} relate that Mike Saks and Peter Winkler observed that known values of $\pi(\PP)$ are all equal to $e(\PP)$. The expression $e(\PP)$ denotes the maximum $m$ such that, for all $n$, the union of the $m$ middle layers in the $n$-dimensional Boolean lattice does not contain $\PP$ as a poset.\footnote{The notation $e(\PP)$ has been used in several places in extremal poset theory. It is not to be confused with the notation for the number of edges in graph $G$, denoted $e(G)$, which we will use later in the paper.}
Bukh~\cite{Bu:2009} verified both Conjecture~\ref{conj:GLL} as well as the Saks-Winkler observation that $\pi(\T)=e(\T)$ for any tree poset $\T$ (meaning that the Hasse diagram of $\T$ is a tree).  Griggs, Li and Lu~\cite{GrLiLu:2012} further verified the stronger conjecture in the case where $\PP$ is a \textit{$k$-diamond poset} for infinitely many values of $k$.  The $k$-diamond poset consists of $k+2$ distinct sets with the property $A\subset B_1,\ldots,B_k\subset C$.


However, a value of $k$ for which they do not determine the $\pi$ value exactly is the case of $k=2$, which we just call the \textit{diamond}.  This poset can also be viewed as the $2$-dimensional Boolean lattice, $\booltwo$.  Hence, the diamond is a poset with four distinct sets $A,B,C,D$ such that $A\subset B,C\subset D$. Griggs, Li and Lu were, however, able to prove that, if the limit exists, then $\pi(\booltwo)\leq 2\frac{3}{11}$.

In Theorem~\ref{thm:main}, we prove that, if the limit exists, then $\pi(\booltwo)\leq 2\frac{1}{4}$.  In other words,
\begin{thm}\label{thm:main}
   Let $\F$ be a diamond-free poset in the $n$-dimensional boolean lattice, $\booln$. Then,
   $$ |\F|\leq (2.25+o(1))\nntwo . $$
\end{thm}

The proof of one of our lemmas (Lemma~\ref{lem:flags}, stated below) was found using Razborov's flag algebra method~\cite{Ra:2007}. After finding an upper bound that can be expressed as a graph invariant, we use the method to show that this invariant is bounded asymptotically by $2.25$. Knowing flag algebras is important to understanding the origins and motivations behind the proof, but it is not necessary for the reader to understand flag algebras in order to understand our proof. In fact, we do not use flag algebra terminology, \textit{per se}.

\section{Proof of Theorem~\ref{thm:main}}

The proof of Theorem~\ref{thm:main} follows from three lemmas: Lemma~\ref{lem:sizelu} (proven in Section~\ref{sec:sizelu}), Lemma~\ref{lem:graph} (proven in Section~\ref{sec:graph}) and Lemma~\ref{lem:flags} (proven in Section~\ref{sec:flags}).

We first use a variant of the YBLM argument that proves Sperner's theorem.  So, we define the Lubell function~\cite{GrLiLu:2012} of a family in a Boolean lattice:
\begin{defn}
If $\F$ is a family of sets in the $n$-dimensional Boolean lattice, the \textbf{Lubell function} of that family is defined to be $\Lam(n,\F)\stackrel{\rmdef}{=}\sum_{F\in\F}\binom{n}{|F|}^{-1}$.

Let $\Lamstar(n,\PP)$ be the maximum of $\Lam(n,\F)$ over all families $\F$ that are both $\PP$-free and contain the empty set. Furthermore, set
$$ \Lamstar(\PP)\stackrel{\rmdef}{=}\limsup_{n\rightarrow\infty}\; \{\Lamstar(n,\PP)\} . $$
\end{defn}~\\

The Lubell function of a family $\mathcal{F}$ gives the average number of times a random full chain will meet $\mathcal{F}$. This will yield an upper bound for the size of a family.
Our main result is that the Lubell function of a $\booltwo$-free poset which contains the empty set is at most $2.25+o(1)$, which immediately implies Theorem~\ref{thm:main} by way of Lemma~\ref{lem:sizelu}. Griggs, Li and Lu~\cite{GrLiLu:2012} address Lemma~\ref{lem:sizelu}.  Here we prove it in detail for completeness.

\begin{lem} \label{lem:sizelu}
Let $\booltwo$ denote the diamond. If $\F$ is $\booltwo$-free in $\booln$ then $|\F|\leq (\Lamstar(\booltwo)+o(1))\nntwo$.  That is, $\La(n,\booltwo)\leq (\Lamstar(\booltwo)+o(1))\nntwo$.
\end{lem}

The expression $\Lamstar(\booltwo)$ provides us with an upper bound on the value of $\La(n,\booltwo)$ but calculating it in general can be difficult. We next introduce a method to compare the value of $\Lamstar(n,\booltwo)$ to a graph invariant that we can calculate more directly. First we introduce some definitions.

\begin{defn}
   For a graph $G$, let $\alpha_i=\alpha_i(G)$ denote the number of three-vertex subgraphs that induce exactly $i$ edges for $i=0,1,2,3$ and let $\beta_j=\beta_j(G)$ denote the number of four-vertex subgraphs that induce exactly $j$ edges for $j=0,\ldots,6$.  If $(X,Y)$ is an ordered bipartition of $V(G)$, then let $\ove(X)$ denote the number of nonedges in the subgraph induced by $X$ and $\ove(Y)$ denote the number of nonedges in the subgraph induced by $Y$.

   Recall that for every nonnegative integer $k$ and real number $n$, $(n)_k=n(n-1)\cdots (n-k+1)$ for $k\geq 1$ and $(n)_0=1$.
\end{defn}

Now we introduce our primary method for bounding the $\Lamstar(n,\F)$ function which will, in turn, bound the size of our extremal function for the size of the family of sets $\F$.
Using the Pochhammer notation, we have the following bound:
\begin{lem}\label{lem:graph}
   For every $\booltwo$-free family $\F$ in $\booln$ with $\emptyset\in\F$, there exist the following:
   \begin{itemize}
      \item a graph $G=(V,E)$ on $v\leq n$ vertices and
      \item a set $W=\{w_{v+1},\ldots,w_n\}$ such that, for each $w\in W$, $(X_w,Y_w)$ is an ordered bipartition of $V$;
   \end{itemize}
   for which
   $$ \Lam(n,\F) \leq 2+ f(n,G,W) , $$
   where, with the notation as above,
   $$ f(n,G,W)=\frac{2\alpha_1(G)-2\alpha_2(G)}{(n)_3}+\frac{6\beta_0(G)}{(n)_4}+\sum_{w\in W}\left[\frac{|X_w|-|Y_w|}{(n)_2}+\frac{4\ove(Y_w)-2\ove(X_w)}{(n)_3}\right] . $$
\end{lem}~\\

Even though the function $f(n,G,W)$ is somewhat complicated, Lemma~\ref{lem:flags} gives that $f$ is bounded asymptotically by $1/4$, regardless of the graph $G$ or the set $W$.
\begin{lem}\label{lem:flags}
   For any integer $n$, graph $G=(V,E)$ on $v\leq n$ vertices and a set $W$, of $n-v$ bipartitions of $V(G)$,
   $$ f(n,G,W)\leq \frac{1}{4}+O\left(\frac{1}{n}\right) . $$
\end{lem}~\\

Theorem~\ref{thm:main} follows because Lemma~\ref{lem:graph} and then Lemma~\ref{lem:flags} give that the Lubell function of any $\booltwo$-free family in $\booln$ that contains the emptyset is at most $2.25+o(1)$.  Lemma~\ref{lem:sizelu} then gives that $|\F|\leq (2.25+o(1))\nntwo$.~\\

We note that Griggs-Li-Lu~\cite{GrLiLu:2012} give two constructions of diamond-free families that both contain the empty set and have Lubell function values of $2+\frac{1}{n(n-1)}\left\lfloor \frac{n^2}{4}\right\rfloor$. Thus it is not possible to obtain a result better than $|\F|\leq \left(2.25+o(1)\right)\nntwo$ by using Lemma~\ref{lem:sizelu}.

\section{Proof of Lemma~\ref{lem:sizelu}}
\label{sec:sizelu}

We first draw upon some conclusions from a paper by Axenovich, Manske and the second author~\cite{AxMaMa:2011}. By Lemma 1 in~\cite{AxMaMa:2011} we know that $\sum_{|k-n/2|\geq n^{2/3}} {n\choose k}\leq 2^{-\Omega(n^{1/3})}\nntwo$. We may assume, therefore that all the elements of $\F$ are close to the middle layer (that is, they are subsets of $[n]$ with size between $n/2-n^{2/3}$ and $n/2+n^{2/3}$). This is because the total number of the remaining elements is $o\big(\nntwo\big)$.

For an event $A$, let $1_A$ be the indicator variable of that event. Let $\mathcal{C}$ be the set of full chains in $\booln$. For a diamond-free family $\F$, let $\F_{\rmmin}$ denote the minimal elements of $\F$.  For any $F'\in\F_{\rmmin}$, let $\mathcal{C}_{F'}$ denote the full chains that contain $F'$ and let $\mathcal{C}_0$ denote the full chains that contain no member of $\F_{\rmmin}$.

Our goal is to bound $\F$ by using the trivial inequality:
\begin{equation}\label{eq:LuBound}
   |\F|\leq\nntwo\sum_{F\in\F}\binom{n}{|F|}^{-1}
\end{equation}
By counting the number of times a full chain contains an element $F$ in $\F$ we see that
\begin{align*}
   \sum_{F\in \F} {n\choose |F|}^{-1} & =
   \frac{1}{n!}\sum_{F\in \F} \sum_{C\in \mathcal{C}}1_{F\in C} \\
   & =  \frac{1}{n!}\sum_{C\in\mathcal{C}_0}\sum_{F\in\F}1_{F\in C}
   +\frac{1}{n!}\sum_{F'\in \F_{\rmmin}}\sum_{C\in \mathcal{C}_{F'}}\sum_{F\in\F}1_{F\in C} .
\end{align*}

Since $\F$ does not have a $4$-chain (this would imply a copy of $\booltwo$), it can be partitioned into at most three antichains: the minimal elements, the maximal elements that are not also minimal and whatever remains.  Hence if $C\in\mathcal{C}_0$, it is the case that $\sum_{F\in\F}1_{F\in\mathcal{C}}\leq 2$ because no chain in $\mathcal{C}_0$ can contain more than two members of $\F$.  So we have the following:
$$ \sum_{F\in \F} {n\choose |F|}^{-1}\leq\frac{2}{n!}|\mathcal{C}_0|
   +\frac{1}{n!}\sum_{F'\in \F_{\rmmin}}\sum_{C\in \mathcal{C}_{F'}}\sum_{F\in\F}1_{F\in C} . $$


One can partition full chains contained in $\mathcal{C}_{F'}$ into collections such that each collection contains all extensions of a given chain from $F'$ to $[n]$. The size of each of these collections is $|F'|!$. So we focus our attention on the interval $\left[F',[n]\right]\stackrel{\rmdef}{=}\left\{S : F'\subseteq S\subseteq [n]\right\}$.  The full chains in this interval are denoted $\mathcal{C}'_{F'}$.  Therefore, $\sum_{C\in \mathcal{C}_{F'}} \sum_{F\in\F} 1_{F\in \mathcal{C}}=|F'|!\sum_{C\in \mathcal{C}'_{F'}}\sum_{F\in \F}1_{F\in\mathcal{C}'_{F'}}$ and so,
\begin{align*}
   \sum_{F\in \F} {n\choose |F|}^{-1} & \leq
   \frac{2}{n!}|\mathcal{C}_0|
   +\frac{1}{n!}\sum_{F'\in \F_{\rmmin}}\frac{(n-|F'|)!|F'|!}{(n-|F'|)!}\sum_{C\in \mathcal{C}'_{F'}}\sum_{F\in\F}1_{F\in C} \\
   & =
   \frac{2}{n!}|\mathcal{C}_0|
   +\sum_{F'\in \F_{\rmmin}}\binom{n}{|F'|}^{-1}\left(\frac{1}{(n-|F'|)!}\sum_{C\in \mathcal{C}'_{F'}}\sum_{F\in\F}1_{F\in C}\right) \\
   & \leq
   \frac{2}{n!}|\mathcal{C}_0|
   +\sum_{F'\in \F_{\rmmin}}\binom{n}{|F'|}^{-1}\Lamstar(n-|F'|,\booltwo) ,
\end{align*}
because $F'$ may be regarded as the empty set in the interval $\left[F',[n]\right]$.

Next we make the observation that $\F_{\rmmin}$ is an antichain and, as such, obeys YBLM.  Set $M_n\stackrel{\rmdef}{=}\max_k\left\{\Lamstar(n-k,\booltwo):|k-n/2|<n^{2/3}\right\}$
\begin{align*}
   \sum_{F\in \F} {n\choose |F|}^{-1} & \leq
   \frac{2}{n!}|\mathcal{C}_0| +\sum_{F'\in \F_{\rmmin}}\binom{n}{|F'|}^{-1}M_n \\
   & \leq \frac{2}{n!}|\mathcal{C}_0| +\left(1-\frac{|\mathcal{C}_0|}{n!}\right)M_n \\
\end{align*}

It is obvious that $\Lamstar(k,\booltwo)$ is at least 2 for $k\geq 2$ and so
\begin{align*}
   \sum_{F\in \F} {n\choose |F|}^{-1} & \leq  M_n=\max_k\left\{\Lamstar(n-k,\booltwo):|k-n/2|<n^{2/3}\right\} \\
   & \leq  \Lamstar(\booltwo)+o(1) .
\end{align*}
By returning to (\ref{eq:LuBound}), we see that $|\F|\leq \left(\Lamstar(\booltwo)+o(1)\right)\nntwo$. This concludes the proof of Lemma~\ref{lem:sizelu}.\hfill~$\Box$~\\

\section{Proof of Lemma \ref{lem:graph}}
\label{sec:graph}

Let us define a partition of the all the chains in $\booln$ as follows. For a family of sets $\F\subseteq\booln$, denote $\Psi_i=\Psi_i(\F)$ to be the set of full chains that contain exactly $i$ members of $\F$ for $i=0,1,\ldots,n+1$. Note that $|\Psi_0|=0$ since $\emptyset \in \F$ and $|\Psi_i|=0$ for $i\geq 4$ since $F$ has no 4-chain. Hence $i\in \{1,2,3\}$.

Since every chain can only hit 1, 2 or 3 elements in a $\booltwo$-free poset $\F$ we have that $n!= |\Psi_3|+|\Psi_2|+|\Psi_1|$. Then for each $i\in \{1,2,3\}$ the total number of times all the chains in $\Psi_i$ contains an element of $\F$ is $i|\Psi_i|$. Recall that $\Lam(n,\F)$ is the average number of times a chain contains an element of $\F$, and hence we have that $\Lam(n,\F)=  \frac{3|\Psi_3|+2|\Psi_2|+|\Psi_1|}{n!}$. Therefore $\Lam(n,\F)= 2+ \frac{|\Psi_3|-|\Psi_1|}{n!}$.

Let $W=\left\{w\in [n]:\{w\}\in \F\right\}$; i.e., the set of singletons in $\F$. Let $V=[n]-W$ and without loss of generality let $W=\{w_{v+1},w_{v+2},\ldots,w_n\}$, where $v=|V|$. Define the graph $G=(V,E(G))$, where
$$ E(G)=\{\{v,v'\} : v,v'\in V\mbox{ and }\{v,v'\}\in\F\}. $$
For every $w\in W$, there is an ordered partition of $V$ called $(X_w,Y_w)$ where
$$ X_w=\left\{x\in V:\{x,w\}\in\F\right\}\qquad\mbox{ and }\qquad Y_w=V-X_w . $$
Observe that if $w,w'$ are distinct members of $W$, then $\{w,w'\}\notin \F$, else we would have a $\booltwo$ with $\emptyset, \{w\},\{w'\}, \{w,w'\}$.~\\

We proceed by placing bounds on both $|\Psi_1|$ and $|\Psi_3|$. First, we find a lower bound on $|\Psi_1|$. We have the following 4 members of $\Psi_1$:
\begin{enumerate}[~~~~(a)]
   \item Let $z_1,z_2,z'\in V$ be a set of 3 vertices that induces exactly two edges with $z_1z_2$ being the nonedge in $G$. Then $\Psi_1$ includes chains  of the form $\emptyset, \{z_1\}, \{z_1,z_2\}, \{z_1,z_2,z'\}, \ldots, [n]$ and chains of the form $\emptyset, \{z_2\}, \{z_2,z_1\}, \linebreak[4] \{z_2,z_1,z'\}, \ldots, [n]$. There are $2\alpha_2(n-3)!$ such chains.  \label{it:alpha2}
   \item Let $w\in W$, $y\in Y_w$ and $z'\in V-\{y\}$ such that $yz'$ is an edge in $G$. Then $\Psi_1$ includes chains of the form $\emptyset,\{y\},\{y,w\},\{y,w,z'\},\ldots,[n]$. There are $\sum_{w\in W}\sum_{y\in Y_w}\deg(y)(n-3)!$ such chains, where $\deg(y)$ is the degree of $y$ in the graph induced by $V$.  \label{it:ywz}
   \item Let $w\in W$, $x\in X_w$ and $z'\in V-\{x\}$ such that $xz'$ is a nonedge in $G$. Then $\Psi_1$ includes chains of the form $\emptyset,\{x\},\{x,z'\},\{x,z',w\},\ldots,[n]$. There are $\sum_{w\in W}\sum_{x\in X_w}\overline{\deg}(x)(n-3)!$ such chains, where $\overline{\deg}(x)$ is the degree of $x$ in the complement of the graph induced by $V$ (i.e., the nondegree of $x$). \label{it:xzw}
   \item Let $w\in W$, $w'\in W-\{w\}$ and $y\in Y_w$. Then $\Psi_1$ includes chains of the form $\emptyset,\{y\},\{y,w\},\{y,w,w'\},\ldots,[n].$ There are $\sum_{w\in W}|Y_w|(|W|-1)(n-3)!$ such chains.\label{it:yww}
\end{enumerate}

It is clear that there is no chain counted twice among the chains (\ref{it:alpha2}), (\ref{it:ywz}), (\ref{it:xzw}) and (\ref{it:yww}). Therefore,
$$ \frac{|\Psi_1|}{(n-3)!}  \geq  2\alpha_2 +\sum_{w\in W}\sum_{y\in Y_w}\deg(y) +\sum_{w\in W}\sum_{x\in X_w}\overline{\deg}(x) +\sum_{w\in W}|Y_w|(|W|-1) .$$

Let us simplify $\sum_{y\in Y_w}\deg(y) +\sum_{x\in X_w}\overline{\deg}(x).$ Fix $w\in W$. Then each edge in $Y_w$ gets counted twice. Each non-edge in $Y_w$ does not get counted. Each edge in $X_w$ does not get counted. The non-edges in $X_w$ each get counted twice. Each edge from $X_w$ to $Y_w$ gets counted once since there is a $y\in Y_w$ incident to the edge. Each non-edge from $Y_w$ to $X_w$ gets counted once since there is a $x\in X_w$ incident to each non-edge. Then we have $\sum_{y\in Y_w}\deg(y) +\sum_{x\in X_w}\overline{\deg}(x)= |X_w||Y_w|+2e(Y_w)+2\ove(X_w)$. Hence,
\begin{equation}\label{in:psione}
   \frac{|\Psi_1|}{(n-3)!} \geq 2\alpha_2 +\sum_{w\in W}\left[|X_w||Y_w|+2e(Y_w)+2\ove(X_w)+(|W|-1)|Y_w|\right].
\end{equation}~\\

Now we turn our attention to $\Psi_3$. It is important to distinguish two types of members of $\F$: We will count separately the types of chains that contain a singleton (that is, chains of the form $\emptyset,\{w\},\ldots,[n]$ for some $w\in W$) and those that do not contain a singleton. Let $\T$ be the minimal nonempty elements of $\F$ and $\U$ be $\F-\{\emptyset\}-\T$. We then have the following 4 possibilities for members of $\Psi_3$.
\begin{enumerate}[~~~~(i)]
   \item Let $w\in W$ and $x\in X_w$. Then $\Psi_3$ contains chains of the form $\emptyset, \{w\}, \{x,w\}, \ldots, [n]$.  We have a total of $\sum_{w\in W}|X_w|(n-3)!$ of this type of chain.\label{it:xw}
   \item Let $w\in W$, $U\in\U$ and $|U|\geq 3$ such that $w\in U$. Then $\Psi_3$ contains chains of the form $\emptyset,\{w\},\ldots,U,\ldots,[n]$.  There are at most $\sum_{w\in W}\left(|Y_w|^2-|Y_w|-2e(Y_w)\right)(n-3)!$ of these types of chains. This bound is proven in Claim~\ref{cl:wUbounds}.\label{it:wU}
   \item Let $v_1,v_2,v'\in V$ be a set of 3 vertices that induces exactly one edge with $v_1v_2$ being the edge in $G$. Then $\Psi_3$ contains chains of the form $\emptyset, \{v_1\}, \{v_1,v_2\}, \{v_1,v_2,v'\}, \ldots, [n]$ and chains of the form $\emptyset, \{v_2\}, \linebreak[4]\{v_2,v_1\}, \{v_2,v_1,v'\}, \ldots, [n]$. The number of chains of this type is $2\alpha_1(n-3)!$. \label{it:alpha1}
   \item Let $U\subset V$ be a member of $\U$ with $|U|\geq 4$ where $U$ does not contain an edge. Then $\Psi_3$ contains chains of the form $\emptyset,\ldots,U,\ldots,[n]$. Note that only one additional member of $\F$ can be in the interval $[\emptyset, U]$, hence at most a $1/|U|$ fraction of these chains contain three members of $\F$. There are at most $\frac{6}{n-3}\betazero (n-3)!$ such chains. This bound will be proven in Claim~\ref{cl:betazero}.\label{it:beta0}
\end{enumerate}~\\

\begin{cl}\label{cl:chainspsithree}
   The chains found in (\ref{it:xw}), (\ref{it:wU}), (\ref{it:alpha1}), and (\ref{it:beta0}) is an exhaustive list of the chains in $\Psi_3$.
\end{cl}

\begin{pfct}{Claim~\ref{cl:chainspsithree}}
Suppose not. Let $C$ be a chain containing three elements of $\F$. Every chain contains $\emptyset$ and in order to have 3 elements of $\F$ in a full chain, that full chain contains an element $T\in \T$. Suppose $T=\{w\}$ for $w\in W$. All such chains are counted in (\ref{it:xw}) and (\ref{it:wU}). Suppose $|T|\geq 2$. Then there must be a $U\in \U$ above $T$ in $C$. If $|T|=2$ then $T$ induces an edge in the graph $G$ and we have a chain of type (\ref{it:alpha1}). If $|T|>2$ then the chain gets counted as a chain of type (\ref{it:beta0}). Therefore, we have counted all possible chains that contain 3 elements of $\F$.
\end{pfct}

\begin{cl}\label{cl:wUbounds}
   For a fixed $w\in W$ the number of chains of type (\ref{it:wU}) is bounded above by $\left(|Y_w|^2-|Y_w|-2e(Y_w)\right)(n-3)!$. Hence, the total number of such chains is at most $\sum_{w\in W}\left(|Y_w|^2-|Y_w|-2e(Y_w)\right)(n-3)!$.
\end{cl}

\begin{pfct}{Claim~\ref{cl:wUbounds}}
For fixed $w\in W$, we can bound the number of chains of type (\ref{it:wU}) as follows:
\begin{align*}
   \lefteqn{\sum_{U\in\U: U\ni w, |U|\geq 3}(|U|-1)!(n-|U|)!} \\
   & =  (n-3)!\sum_{U\in\U: U\ni w, |U|\geq 3}\binom{|Y_w|}{|U|-1}^{-1}\frac{(|Y_w|)_{|U|-1}}{(n-3)_{|U|-3}} \\
   & =  (n-3)!|Y_w|(|Y_w|-1)\sum_{U\in\U: U\ni w, |U|\geq 3}\binom{|Y_w|}{|U|-1}^{-1}\frac{(|Y_w|-2)_{|U|-3}}{(n-3)_{|U|-3}}
\end{align*}

Since $|Y_w|\leq n-1$, we have $\frac{(|Y_w|-2)_{|U|-3}}{(n-3)_{|U|-3}}\leq 1$. Furthermore, consider subsets of $Y_w$. Note that $U-\{w\}\subseteq Y_w$, else we have a chain of length 4 in $\F$.  The family $\{U-\{w\} : U\in\U, U\ni w, |U|\geq 3\}\cup E\left(\left.G\right|_{Y_w}\right)$ forms an antichain in $\bool_{|Y_w|}$. First, no edge $y_1y_2$ can be a subset of one of those sets $U-\{w\}$, otherwise $\emptyset\subset \{w\},\{y_1,y_2\}\subset U$. Second, any two sets $U_1\subset U_2$ in this family would form a chain of length 4 in $\F$, namely, $\emptyset\subset\{w\}\subset U_1\subset U_2$. So, the YBLM inequality gives $e(Y_w)\binom{|Y_w|}{2}^{-1}+\sum_{U\in\U: U\ni w, |U|\geq 3}\binom{|Y_w|}{|U|-1}^{-1}\leq 1$. Thus,
\begin{align*}
   \lefteqn{\sum_{U\in\U: U\ni w, |U|\geq 3}(|U|-1)!\cdot |U|!(n-|U|)!} \\
   & \leq (n-3)!|Y_w|(|Y_w|-1)\sum_{U\in\U: U\ni w, |U|\geq 3}\binom{|Y_w|}{|U|-1}^{-1} \\
   & \leq  (n-3)!|Y_w|(|Y_w|-1)\left(1-e(Y_w)\binom{|Y_w|}{2}^{-1}\right) \\
   & = (n-3)!\left(|Y_w|^2-|Y_w|-2e(Y_w)\right)
\end{align*}
This concludes the proof of Claim~\ref{cl:wUbounds}.
\end{pfct}

\begin{cl}\label{cl:betazero}
   The number of chains of type (\ref{it:beta0}) is at most $\frac{6}{n-3}\betazero (n-3)!$.
\end{cl}

\begin{pfct}{Claim~\ref{cl:betazero}}
If we fix some $T\in\T$ with $|T|\geq 3$, then we will count these chains by replacing the members of $\U$ that are supersets of $T\in\T$ with supersets of $T$ that cover $T$. A set $U$ \textit{covers} set $T$ if there is a $u\in U$ such that $T=U-\{u\}$. This approach is also discussed in~\cite{AxMaMa:2011}.

Let $\U(T)$ denote the members of $\U$ that are supersets of $T$ and $\U_0(T)$ denote the family of all sets $U$ such that $U$ covers $T$ but there is no other member of $\T$ that $U$ covers.  Note that members of $\U_0(T)$ need not be members of $\U$.  Note further that, by definition, each $\U_0(T)$ is an antichain.

If we consider the family $\U'=\bigcup_{T\in\T}\U_0(T)$, then $\{\emptyset\}\cup\T\cup\U'$ has no $\booltwo$.  This is because, first, no $U\in \U'$ can be a superset of distinct $T_1,T_2\in\T$ by definition.  Second, if $U_1,U_2\in\U'$ such that $U_1\subset U_2$, then there are distinct $T_1,T_2$ such that $T_1\subset U_1$ and $T_2\subset U_2$, but then $U_2\supset T_1$, a contradiction to the definition of $\U_0(T_2)$.

Furthermore, any chain that contains $T$ and some $U\in\U(T)$ has the property that there is some $U_0\in\U_0(T)$ such that this chain also contains $U_0$. If not, then the chain contains some $X$ that covers $T$ for which $T\subset X\subset U$.  The only reason that $X\not\in\U_0(T)$ is that there is another $T'\in\T$ such that $X\supset T'$.  Thus, $\emptyset\subset T,T'\subset U$, a $\booltwo$.

So, we can bound the number of chains of type (\ref{it:beta0}) as follows:
\begin{align*}
   \lefteqn{\hspace{-.5in} \sum_{U\in\U, U\subseteq V, |U|\geq 4} \frac{1}{|U|}|U|!(n-|U|)!} \\
   \hspace{.5in} & =  (n-3)!\sum_{U\in\U, U\subseteq V, |U|\geq 4} \frac{1}{|U|}\cdot\frac{|U|!}{(n-3)_{|U|-3}} \\
   & = (n-3)!\frac{(v)_4}{(n-3)}\sum_{U\in\U, U\subseteq V, |U|\geq 4} \frac{1}{|U|}\binom{v}{|U|}^{-1}\frac{(v-4)_{|U|-4}}{(n-4)_{|U|-4}} .
\end{align*}

We have that $v\leq n$, hence $\frac{(v-4)_{|U|-4}}{(n-4)_{|U|-4}}\leq 1$. Furthermore, $1/|U|\leq 1/4$ for all $U$ in the sum.  Therefore,
$$ \sum_{U\in\U, U\subseteq V, |U|\geq 4} \frac{1}{|U|}|U|!(n-|U|)!
   \leq (n-3)!\frac{(v)_4}{(n-3)}\sum_{U\in\U, U\subseteq V, |U|\geq 4}\frac{1}{4}\binom{v}{|U|}^{-1} . $$

The members of $\{U\in\U : U\subseteq V, |U|\geq 4\}$ together with quadruples of vertices that induce at least one edge form an antichain.  Hence, the YBLM inequality gives that $\left(\binom{v}{4}-\beta_0\right)\binom{v}{4}^{-1}+\sum\binom{v}{|U|}^{-1}\leq 1$.  Hence,
\begin{align*}
   \sum_{U\in\U, U\subseteq V, |U|\geq 4} \frac{1}{|U|}|U|!(n-|U|)!
   & \leq  (n-3)!\frac{(v)_4}{(n-3)}\sum_{U\in\U, U\subseteq V, |U|\geq 4}\frac{1}{4}\binom{v}{|U|}^{-1} \\
   & \leq  (n-3)!\frac{(v)_4}{4(n-3)}\betazero\binom{v}{4}^{-1}  \\
   & =  (n-3)!\frac{6}{n-3}\betazero .
\end{align*}

This concludes the proof of Claim~\ref{cl:betazero}.
\end{pfct}~\\

Returning to the proof of Lemma~\ref{lem:graph}, we can combine the bounds on $\Psi_1$ and $\Psi_3$:
\begin{align*}
   \frac{|\Psi_3|-|\Psi_1|}{(n-3)!} & \leq  2\alpha_1+\frac{6}{n-3}\beta_0+\sum_{w\in W}\left[|X_w|(n-2)+|Y_w|(|Y_w|-1)-2e(Y_w)\right] \\
   &  -2\alpha_2-\sum_{w\in W}\left[|X_w||Y_w|+2e(Y_w)+2\ove(X_w)+(|W|-1)|Y_w|\right] \\
   & =   2(\alpha_1-\alpha_2)+\frac{6}{n-3}\beta_0 \\
   & \hspace{.45in} +\sum_{w\in W}\left[|X_w|(n-2) +(|Y_w|)_2 -|Y_w|(|W|-1) -|X_w||Y_w|\right. \\
   &  \hspace{.85in}\left.-4e(Y_w) -2\ove(X_w)\right] \\
   & =   2(\alpha_1-\alpha_2)+\frac{6}{n-3}\beta_0 \\
   & \hspace{.45in} +\sum_{w\in W}\left[|X_w|(n-2) +|Y_w|(|Y_w|-|W|) -|X_w||Y_w|\right. \\
   &  \hspace{.85in}\left.-2(|Y_w|)_2+4\ove(Y_w)-2\ove(X_w)\right]
\end{align*}

We can collect terms and observe that $|Y_w|\left(|Y_w|-|W|-|X_w|-2|Y_w|+2\right)=-|Y_w|(n-2)$. The expression then greatly simplifies to give us the following:
\begin{align*}
   \frac{|\Psi_3|-|\Psi_1|}{(n-3)!} & \leq  2(\alpha_1-\alpha_2)+\frac{6}{n-3}\beta_0 \\
   &  \hspace{.25in}+\sum_{w\in W}\left[\left(|X_w|-|Y_w|\right)(n-2) +4\ove(Y_w)-2\ove(X_w)\right] .
\end{align*}

If we divide by $(n)_3$ we get the bound on $3|\Psi_3|+2|\Psi_2|+|\Psi_1|=2+|\Psi_3|-|\Psi_1|$, as needed. This concludes the proof of Lemma~\ref{lem:graph}.\hfill~$\Box$~\\

\section{Proof of Lemma~\ref{lem:flags}}
\label{sec:flags}

\newcommand{\hset}{\mathcal{H}_4}
The first thing we do is compute $f(n,G,W)$ by a summation on the set of induced subgraphs on $4$ vertices.  In order to do so, however, we must eliminate the cases of when $v\leq 3$.

If $v\leq 3$, then $\alpha_1\leq 1$, $|X_w|\leq 3$ and $\ove(Y_w)\leq 3$.  So trivially,
$$ f(n,G,W)\leq\frac{2}{(n)_3}+n\left(\frac{3}{(n)_2}+\frac{12}{(n)_3}\right) , $$
which is at most $1/4$ if $n\geq 11$.

Now we assume that $v\geq 4$.  Let $\hset=\binom{V}{4}$, the set of all 4-tuples of vertices.  For ease of notation, if $H\in\hset$, we will use $H$ to mean both that set of 4 vertices and the subgraph induced by them.

We observe that $\sum_{H\in\hset}\alpha_i(H)=(v-3)\alpha_i(G)$ for $i=0,1,2,3$ and for any set $S\subseteq V$, we have $\sum_{H\in\hset}\ove(S\cap H)=\binom{v}{2}\ove(S)$ and $\sum_{H\in\hset}|S\cap H|=\binom{v}{3}|S|$.  So, we can rewrite $f(n,G,W)$ as follows:
\begin{align*}
   f(n,G,W) & =  \frac{2\alpha_1-2\alpha_2}{(n)_3}+\frac{6\beta_0}{(n)_4}+\sum_{w\in W}\left[\frac{|X_w|-|Y_w|}{(n)_2}+\frac{4\ove(Y_w)-2\ove(X_w)}{(n)_3}\right] \\
   & =  \frac{1}{\binom{v}{4}}\sum_{H\in\hset}\left[\frac{(v)_3}{3(n)_3}\cdot\frac{\alpha_1(H)-\alpha_2(H)}{4} +\frac{(v)_4}{4(n)_4}\beta_0(H)\right. \\
   &  \hspace{.65in}\left.+\sum_{w\in W}\left(\frac{v}{(n)_2}\cdot\frac{|X_w\cap H|-|Y_w\cap H|}{4}\right.\right. \\
   &  \hspace{1.15in}\left.\left.+\frac{(v)_2}{2(n)_3}\cdot\frac{4\ove(Y_w\cap H)-2\ove(X_w\cap H)}{6}\right)\right] .
\end{align*}

We will see in Claim~\ref{cl:eps} that in most cases of $H\in\hset$, the largest value of the summand occurs when $X_w\cap H=H$.  So, we will rearrange the terms:
\begin{align}
   f(n,G,W) & =  \frac{1}{\binom{v}{4}}\sum_{H\in\hset}\left[\frac{(v)_3}{3(n)_3}\cdot\frac{\alpha_1(H)-\alpha_2(H)}{4} +\frac{(v)_4}{4(n)_4}\beta_0(H)+\frac{(n-v)v}{(n)_2}\right. \nonumber \\
   &  \hspace{.65in}\left. +\frac{(n-v)(v)_2\,\ove(H)}{6(n)_3} +\sum_{w\in W}\epsilon(n,w,G,H)\right] , \label{eq:fH}
\end{align}
where
$$ \epsilon(n,w,G,H)=-\frac{v}{(n)_2}\cdot\frac{2|Y_w\cap H|}{4} +\frac{(v)_2}{2(n)_3}\cdot\frac{2\ove(H)+4\ove(Y_w\cap H)-2\ove(X_w\cap H)}{6} . $$~\\

At this point, we need notation previously used in ~\cite{AxMaMa:2011} to describe the eleven distinct nonisomorphic graphs on exactly 4 vertices as follows.
\begin{notat}
   For $i=0,1,5,6$, the $4$-vertex graph with exactly $i$ edges is denoted $H_i$.

   The graph with exactly two edges that are incident is $\hwedge$ and the graph with exactly two edges that are nonincident is $\hpar$. Their complements are $\hq$ and $\hcyc$, respectively.

   The graph inducing a star with three edges (the claw) is $\hw$, the graph inducing a triangle is $\htri$ and the path with three edges is $\hpath$.

   We use $\cong$ to denote that two graphs are isomorphic.
\end{notat}~\\

Claim~\ref{cl:eps} bounds the $\epsilon$ term in (\ref{eq:fH}).  Note that it is $0$ for all but three of the eleven $4$-vertex graphs.
\begin{cl}\label{cl:eps}
For any integer $n$, graph $G=(V,E)$ on $v<n$ vertices, ordered bipartition $(X_w,Y_w)$ of $V$ and $H\in\hset$,
$$ \epsilon(n,w,G,H) \leq \begin{cases}
                                    \frac{3v}{(n)_2}\max\{0,\frac{v-1}{n-2}-\frac{2}{3}\} , & \mbox{if $H\cong\hzero$;} \\
                                    \frac{5v}{2(n)_2}\max\{0,\frac{v-1}{n-2}-\frac{4}{5}\} , & \mbox{if $H\cong\hone$;} \\
                                    \frac{5v}{3(n)_2}\max\{0,\frac{v-1}{n-2}-\frac{9}{10}\} , & \mbox{if $H\cong\hwedge$;} \\
                                    0 , & \mbox{otherwise.}
                                 \end{cases} $$
\end{cl}

\begin{pfct}{Claim~\ref{cl:eps}}
We note that we have defined $\epsilon(n,w,G,H)$ to be zero if $Y_w\cap H=\emptyset$.

For each case of $|Y_w\cap H|$, we choose the set of vertices that makes the expression $4\ove(Y_w\cap H)-2\ove(X_w\cap H)$ as large as possible. We look at each individual possibility for $|Y_w\cap H|$.  First suppose $\ove(H)\leq 3$. In this case, $\ove(Y_w\cap H)$ is at most $\ove(H)$ if $|Y_w\cap H|\in\{3,4\}$.  In all such cases, $\epsilon(n,w,G,H)\leq 0$, see Table~\ref{table:eps:3} for the table of terms in the case $\ove(H)\leq 3$.

We can now consider the remaining four possibilities for $H$.  The case $\hpar$ is displayed in Table~\ref{table:eps:hpar}, in which it is shown that $\epsilon(n,w,G,\hpar)\leq 0$.  The case $\hwedge$ is displayed in Table~\ref{table:eps:hwedge}, in which it is shown that $\epsilon(n,w,G,\hpar)\leq 0$ unless $|Y_w\cap H|=3$ and $v$ is large enough.  The case $\hone$ is displayed in Table~\ref{table:eps:hone} in which $\epsilon(n,w,G,\hone)$ is either at most $0$ or at most the case where $|Y_w\cap H|=4$.  The case $\hzero$ is displayed by Table~\ref{table:eps:hzero} in which $\epsilon(n,w,G,\hzero)$ is either at most $0$ or at most the case where $|Y_w\cap H|=4$.

This concludes the proof of Claim~\ref{cl:eps}.
\end{pfct}

Note that for fixed $n$, $G$ and $W$, we have expressed $f(n,G,W)$ as $\frac{1}{\binom{v}{4}}\sum_H d(H)$ for some density function $d$ of $H$ in Equation (\ref{eq:fH}).  Moreover, we can use the bounds for the $\epsilon$ terms in Claim~\ref{cl:eps} in order to get upper bounds for $d(H)$, call them $d^*(H)$.  As such, we have Table~\ref{table:dstar}, which details the values of $d^*(H)$ for the eleven different nonisomorphic graphs on $4$ vertices.~\\

\noindent\textbf{Case 1.} $4\leq v\leq (2n-1)/3$.~\\

In this case, the expressions involving the $\max$ function are zero.  These simplified values of $d^*(H)$ are displayed in Table~\ref{table:dplusc_one}.

It is trivial to see that, if $v\leq (2n-1)/3$ and $H\in\{\hw,\hpath,\hcyc,\hq,\hfive,\hsix\}$ that $d^*(H)\leq\frac{(n-v)v}{(n)_2}$.  For the rest of the graphs $H$, the expressions in Table~\ref{table:dplusc_one} can be simplified as follows:
\begin{align*}
   d^*(\hzero) & \leq  \ts\frac{(n-v)v}{(n)_2} +\frac{(v)_2}{4(n)_4}\left(v^2+4nv-17v-4n^2+12n+6\right) \\
   d^*(\hone) & \leq  \ts\frac{(n-v)v}{(n)_2} +\frac{(v)_2}{6(n)_3}(6v-5n-2) \\
   d^*(\hwedge) & \leq  \ts\frac{(n-v)v}{(n)_2} +\frac{(v)_2}{12(n)_3}(9v-8n-2) \\
   d^*(\hpar) & \leq \ts\frac{(n-v)v}{(n)_2} +\frac{(v)_2}{3(n)_3}(3v-2n-2) \\
   d^*(\htri) & \leq \ts\frac{(n-v)v}{(n)_2} +\frac{(v)_2}{4(n)_3}(3v-2n-2) \\
\end{align*}

With this simplification, it is easy to see that $d^*(H)<\frac{(n-v)v}{(n)_2}$ if $v\leq (2n-1)/3$ and $H\in\{\hone,\hwedge,\hpar,\htri\}$.  As to $\hzero$, the expression $v^2+4nv-17v-4n^2+12n+6$ increases in $v$ for all $n\geq 4$.  Hence, its maximum value for $v\leq (2n-1)/3$ is at most $(-8n^2-10n+106)/9$, which is negative for all $n\geq 4$.

Therefore, if $4\leq v\leq (2n-1)/3$ and $n\geq 4$, then
$$ f(n,G,W)\leq \max_{H\in\hset}\{d^*(H)\}\leq \frac{(n-v)v}{(n)_2}\leq\frac{1}{(n)_2}\left\lfloor\frac{n^2}{4}\right\rfloor . $$~\\

\noindent\textbf{Case 2.} $v\geq 2n/3$ (and $v\geq 4$).~\\

We will upper bound $f(n,G,W)$ by adding a nonnegative term to $\binom{v}{4}^{-1}\sum d^*(H)$. We use $N(v)$ to denote the neighborhood of vertex $v$, $N[v]=N(v)\cup\{v\}$ to denote the closed neighborhood of vertex $v$. For any vertex-subset $S$, we denote $\overline{S}=V-S$.  Now we can bound $f(n,G,W)$, where the coefficient $\gamma\geq 0$ will be chosen later:
\begin{align}
   f(n,G,W) & \leq  \frac{1}{\binom{v}{4}}\sum_{H\in\hset}d^*(H) \nonumber \\
   & \leq  \frac{1}{\binom{v}{4}}\sum_{H\in\hset}d^*(H) \nonumber \\
   &  \quad+\frac{\gamma}{\binom{v}{4}}\sum_{(z_1,z_2) : z_1z_2\in E(\overline{G})}\left(|N(z_1)\cap\overline{N[z_2]}| -|\overline{N[z_1]}\cap N(z_2)|\right)^2 \label{eq:sq:1} \\
   & \quad +\frac{\gamma}{\binom{v}{4}}\sum_{(z_1,z_2) : z_1z_2\in E(G)}\left(|N(z_1)\cap N(z_2)| -|\overline{N[z_1]}\cap \overline{N[z_2]}|\right)^2 \label{eq:sq:2}
\end{align}

So, the expression (\ref{eq:sq:1}) is a sum over ordered pairs that form nonedges in $G$ and so $z_1\neq z_2$ in this summation.  The expression (\ref{eq:sq:2}) is a sum over ordered pairs that form edges in $G$.
The reader familiar with Razborov's flag algebra method~\cite{Ra:2007} will see it at work here.  There are two ``types'' -- the edges and nonedges.  The flags are subgraphs induced by three vertices.  For each type, the positive semidefinite matrix is a $4\times 4$ matrix that is a multiple of a matrix with two diagonal entries equal to $1$, two off-diagonal entries equal to $-1$ and the remaining twelve entries being equal to $0$.  For more on flag algebras, we refer the reader to Baber and Talbot~\cite{BaTa:2011} who give a nice discussion of an application to hypergraphs.

The key observation regarding expressions (\ref{eq:sq:1}) and (\ref{eq:sq:2}) is that these terms can be rewritten so that they count an expression we call $\gamma c(H)$ for each subgraph $H$ on $4$ vertices, as well as some terms corresponding to subgraphs on $3$ vertices.  For example, the expression (\ref{eq:sq:1}) assigns $-16\gamma$ for each subgraph isomorphic to $\hpar$, because it assigns a value of $-2\gamma$ for each ordered pair of distinct nonadjacent vertices, of which there are eight.  The expression (\ref{eq:sq:2}) assigns $+8\gamma$ for each subgraph isomorphic to $\hpar$, because it assigns a value of $+2\gamma$ for each ordered pair of distinct adjacent vertices, of which there are four.  This gives a net value of $-8\gamma$.  We compute all values of $\gamma c(H)$ in this way and display them in Table~\ref{table:see}.  This gives the expression:
\begin{align}
   f(n,G,W) &\leq  \frac{1}{\binom{v}{4}}\sum_{H\in\hset}d^*(H) +\frac{\gamma}{\binom{v}{4}}\sum_{H\in\hset}c(H) +\frac{\gamma}{\binom{v}{4}}(6\alpha_1(G)+6\alpha_3(G)) \\
   & \leq  \frac{1}{\binom{v}{4}}\sum_{H\in\hset}\left(d^*(H)+\gamma c(H)\right) +\frac{24\gamma}{v-3} . \label{eq:bound}
\end{align}

We choose the value of $\gamma$ so as to ensure that $d^*(\hsix)+\gamma c(\hsix)=\frac{(n-v)v}{(n)_2}+24\gamma$ is exactly $1/4$.  Hence, $\gamma=\frac{1}{96}-\frac{(n-v)v}{24(n)_2}$.  Since $v\geq 2n/3$, this choice of $\gamma$ is nonnegative as long as $n\geq 7$. It remains to show that $d^*(H)+\gamma c(H)\leq 1/4$ for the remaining choices of $H$.

In Table~\ref{table:dplusc_twoa}, we list all of the expressions for $d^*(H)+\gamma c(H)$ and simplify it to obtain Table~\ref{table:dplusc_twob}, which gives $d^*(H)+\gamma\cdot c(H)-1/4$ for each of the eleven graphs $H$ on $4$ vertices.  We now only need to see that these expressions are at most $0$.

First, we see in Table~\ref{table:dplusc_twob} that $d^*(H_0)+\gamma c(H_0)-1/4\leq 0$.  As to the rest, the easiest approach is to express the functions in Table~\ref{table:dplusc_twob} in terms of $x=v/n$, ignoring the lower-order terms and obtaining a function $g_H(x)$. Since the numerators are polynomials in $v$ and $n$ and denominators are polynomials in $n$, then $g_H(x)\leq 0$ implies that $d^*(H_0)+\gamma c(H_0)-1/4=O(1/n)$.

The functions $g_H(x)$ are given in Table~\ref{table:xfunct} along with their maximum value over the range $2/3\leq x\leq 1$.  They are strictly negative for all $H\not\in\{\hpar,\htri\}$ and at most zero for those. This concludes the proof of Lemma~\ref{lem:flags}.\hfill~$\Box$~\\

\section{Concluding remarks}
Griggs, Li and Lu conjectured that the Lubell function for a diamond-free family in $\booln$ is at most $2+\frac{1}{(n)_2}\left\lfloor\frac{n^2}{4}\right\rfloor$.  For $n$ large enough, we have shown this in the case where $v<\frac{2n}{3}$.

In the case of $v\geq\frac{2n}{3}$, for all $H\not\in\{\hpar,\htri\}$, the fact that $g_H(x)$ is always strictly negative (Table~\ref{table:xfunct}) gives that $d^*(H)+\gamma c(H)-1/4\leq 0$ for $n$ large enough.
Furthermore, it can be shown that $d^*(\hpar)+\gamma c(\hpar)-1/4\leq 0$ and $d^*(\htri)+\gamma c(\htri)-1/4\leq 0$ with equality if and only if $v=n$.

Coupling this with the fact that $\frac{24\gamma}{v-3}=\frac{24}{v-3}\left(\frac{1}{96}-\frac{(n-v)v}{24(n)_2}\right)\leq\frac{1}{4(n-3)}$, the $O(1/n)$ term in the statement of the lemma can be seen to be as small as $1/(4n-12)$ for $n$ large enough.

However,   $\frac{1}{(n)_2}\left\lfloor\frac{n^2}{4}\right\rfloor-\frac{1}{4}$ is either $\frac{1}{4n}$ or $\frac{1}{4(n-1)}$, depending on the parity of $n$.  We believe that the bound by Griggs, Li and Lu is correct. That is, if $\F$ is a diamond free family, then $\Lambda(n,\F)\leq 2+ \frac{1}{n(n-1)}\left\lfloor\frac{n^2}{4}\right\rfloor$.~\\

\section{Acknowledgements}
The authors would like to thank Peter Keevash and John Talbot for suggesting the method of flag algebras.~\\

\section*{References}
\bibliographystyle{plain}
\bibliography{bibliographyfile}

\appendix
\FloatBarrier
\section{Tables}

\subsection{Tables of values to prove Claim~\ref{cl:eps}.}
\FloatBarrier
\begin{table}[ht]\centering
\begin{tabular}{||r|rclr||}\hline
$|Y_w\cap H|$      & \multicolumn{4}{l||}{maximum value of $\epsilon(n,w,G,H)$ if $\ove(H)\leq 3$.} \\ \hline\hline
1 & $\ts -\frac{v}{2(n)_2}+\frac{(v)_2}{2(n)_3}\frac{2\ove(H)}{6}$              & $\ts \leq$ & $\ts \frac{v}{2(n)_2}\left(-1+\frac{v-1}{n-2}\cdot\frac{\ove(H)}{3}\right)$ & $\ts \leq 0$ \\ \hline
2 & $\ts -\frac{v}{(n)_2}+\frac{(v)_2}{2(n)_3}\frac{2\ove(H)+4\cdot 1}{6}$      & $\ts \leq$ & $\ts \frac{v}{(n)_2}\left(-1+\frac{v-1}{n-2}\cdot\frac{\ove(H)+2}{6}\right)$ & $\ts <0$ \\ \hline
3 & $\ts -\frac{3v}{2(n)_2}+\frac{(v)_2}{2(n)_3}\frac{2\ove(H)+4\ove(H)}{6}$    & $\ts =$    & $\ts \frac{v}{(n)_2}\left(-\frac{3}{2}+\frac{v-1}{n-2}\cdot\frac{\ove(H)}{2}\right)$ & $\ts \leq 0$ \\ \hline
4 & $\ts -\frac{2v}{(n)_2}+\frac{(v)_2}{2(n)_3}\frac{2\ove(H)+4\ove(H)}{6}$     & $\ts =$    & $\ts \frac{v}{(n)_2}\left(-2+\frac{v-1}{n-2}\cdot\frac{\ove(H)}{2}\right)$ & $\ts <0$ \\ \hline
\end{tabular}
\caption{Maximum value of $\epsilon(n,w,G,H)$ if $\ove(H)\leq 3$.}
\label{table:eps:3}
\end{table}

\begin{table}[ht]\centering
\begin{tabular}{||r|rclr||}\hline
$|Y_w\cap H|$      & \multicolumn{4}{l||}{maximum value of $\epsilon(n,w,G,\hpar)$} \\ \hline\hline
1 & $\ts -\frac{v}{2(n)_2}+\frac{(v)_2}{2(n)_3}\frac{2\cdot 4-2\cdot 2}{6}$         & $\ts =$ & $\ts \frac{v}{2(n)_2}\left(-\frac{1}{2}+\frac{1}{3}\frac{v-1}{n-2}\right)$ & $\ts <0$ \\ \hline
2 & $\ts -\frac{v}{(n)_2}+\frac{(v)_2}{2(n)_3}\frac{2\cdot 4+4\cdot 1-2\cdot 1}{6}$ & $\ts =$ & $\ts \frac{v}{(n)_2}\left(-1+\frac{5}{6}\frac{v-1}{n-2}\right)$ & $\ts <0$ \\ \hline
3 & $\ts -\frac{3v}{2(n)_2}+\frac{(v)_2}{2(n)_3}\frac{2\cdot 4+4\cdot 3}{6}$        & $\ts =$ & $\ts \frac{v}{(n)_2}\left(-\frac{3}{2}+\frac{5}{3}\frac{v-1}{n-2}\right)$ & $\ts <0$ \\ \hline
4 & $\ts -\frac{2v}{(n)_2}+\frac{(v)_2}{2(n)_3}\frac{6\cdot 4}{6}$                  & $\ts =$ & $\ts \frac{v}{(n)_2}\left(-2+2\frac{v-1}{n-2}\right)$ & $\ts \leq 0$ \\ \hline
\end{tabular}
\caption{Maximum value of $\epsilon(n,w,G,\hpar)$.}
\label{table:eps:hpar}
\end{table}

\begin{table}[ht]\centering
\begin{tabular}{||r|rclr||}\hline
$|Y_w\cap H|$      & \multicolumn{4}{l||}{maximum value of $\epsilon(n,w,G,\hwedge)$} \\ \hline\hline
1 & $\ts -\frac{v}{2(n)_2}+\frac{(v)_2}{2(n)_3}\frac{2\cdot 4-2\cdot 1}{6}$  & $\ts =$ & $\ts \frac{v}{2(n)_2}\left(-\frac{1}{2}+\frac{1}{2}\frac{v-1}{n-2}\right)$ & $\ts \leq 0$ \\ \hline
2 & $\ts -\frac{v}{(n)_2}+\frac{(v)_2}{2(n)_3}\frac{2\cdot 4+4\cdot 1}{6}$   & $\ts =$ & $\ts \frac{v}{(n)_2}\left(-1+\frac{v-1}{n-2}\right)$ & $\ts <0$ \\ \hline
3 & $\ts -\frac{3v}{2(n)_2}+\frac{(v)_2}{2(n)_3}\frac{2\cdot 4+4\cdot 3}{6}$ & $\ts =$ & $\ts \frac{5v}{3(n)_2}\left(-\frac{9}{10}+\frac{v-1}{n-2}\right)$ & \\ \hline
4 & $\ts -\frac{2v}{(n)_2}+\frac{(v)_2}{2(n)_3}\frac{6\cdot 4}{6}$           & $\ts =$ & $\ts \frac{v}{(n)_2}\left(-2+2\frac{v-1}{n-2}\right)$ & $\ts \leq 0$ \\ \hline
\end{tabular}
\caption{Maximum value of $\epsilon(n,w,G,\hwedge)$.}
\label{table:eps:hwedge}
\end{table}

\begin{table}[ht]\centering
\begin{tabular}{||r|rclr||}\hline
$|Y_w\cap H|$      & \multicolumn{4}{l||}{maximum value of $\epsilon(n,w,G,\hone)$} \\ \hline\hline
1 & $\ts -\frac{v}{2(n)_2}+\frac{(v)_2}{2(n)_3}\frac{2\cdot 5-2\cdot 2}{6}$  & $\ts =$ & $\ts \frac{v}{2(n)_2}\left(-\frac{1}{2}+\frac{1}{2}\frac{v-1}{n-2}\right)$ & $\ts \leq 0$ \\ \hline
2 & $\ts -\frac{v}{(n)_2}+\frac{(v)_2}{2(n)_3}\frac{2\cdot 5+4\cdot 1}{6}$   & $\ts =$ & $\ts \frac{7v}{6(n)_2}\left(-\frac{6}{7}+\frac{v-1}{n-2}\right)$ & \\ \hline
3 & $\ts -\frac{3v}{2(n)_2}+\frac{(v)_2}{2(n)_3}\frac{2\cdot 5+4\cdot 3}{6}$ & $\ts =$ & $\ts \frac{11v}{6(n)_2}\left(-\frac{9}{11}+\frac{v-1}{n-2}\right)$ & \\ \hline
4 & $\ts -\frac{2v}{(n)_2}+\frac{(v)_2}{2(n)_3}\frac{6\cdot 5}{6}$           & $\ts =$ & $\ts \frac{5v}{2(n)_2}\left(-\frac{4}{5}+\frac{v-1}{n-2}\right)$ & \\ \hline
\end{tabular}
\caption{Maximum value of $\epsilon(n,w,G,\hone)$.}
\label{table:eps:hone}
\end{table}

\begin{table}[ht]\centering
\begin{tabular}{||r|rclr||}\hline
$|Y_w\cap H|$      & \multicolumn{4}{l||}{maximum value of $\epsilon(n,w,G,\hzero)$} \\ \hline\hline
1 & $\ts -\frac{v}{2(n)_2}+\frac{(v)_2}{2(n)_3}\frac{2\cdot 4-2\cdot 3}{6}$         & $\ts =$ & $\ts \frac{v}{2(n)_2}\left(-\frac{1}{2}+\frac{1}{2}\frac{v-1}{n-2}\right)$ & $\ts \leq 0$ \\ \hline
2 & $\ts -\frac{v}{(n)_2}+\frac{(v)_2}{2(n)_3}\frac{2\cdot 4+4\cdot 1-2\cdot 1}{6}$ & $\ts =$ & $\ts \frac{7v}{6(n)_2}\left(-\frac{6}{7}+\frac{v-1}{n-2}\right)$ & \\ \hline
3 & $\ts -\frac{3v}{2(n)_2}+\frac{(v)_2}{2(n)_3}\frac{2\cdot 4+4\cdot 3}{6}$        & $\ts =$ & $\ts \frac{2v}{(n)_2}\left(-\frac{3}{4}+\frac{v-1}{n-2}\right)$ & \\ \hline
4 & $\ts -\frac{2v}{(n)_2}+\frac{(v)_2}{2(n)_3}\frac{6\cdot 6}{6}$                  & $\ts =$ & $\ts \frac{3v}{(n)_2}\left(-\frac{2}{3}+\frac{v-1}{n-2}\right)$ & \\ \hline
\end{tabular}
\caption{Maximum value of $\epsilon(n,w,G,\hzero)$.}
\label{table:eps:hzero}
\end{table}
\FloatBarrier

\FloatBarrier
\subsection{Tables for the proof of Lemma~\ref{lem:flags}}
\begin{table}[ht]\centering
\begin{tabular}{||r|llll||}\hline
$H$      & \multicolumn{4}{l||}{$d^*(H)$} \\ \hline\hline
$\hzero$ & $\ts\frac{1}{4}\frac{(v)_4}{(n)_4}$   & $\ts +\frac{(n-v)v}{(n)_2}$      & $\ts -\frac{(n-v)(v)_2}{(n)_3}$ & $\ts +\frac{3(n-v)v}{(n)_2}\max\left\{0, \frac{v-1}{n-2}-\frac{2}{3}\right\}$ \\ \hline
$\hone$ & $\ts\frac{1}{6}\frac{(v)_3}{(n)_3}$    & $\ts +\frac{(n-v)v}{(n)_2}$      & $\ts -\frac{5(n-v)(v)_2}{6(n)_3}$ & $\ts +\frac{5(n-v)v}{2(n)_2}\max\left\{0, \frac{v-1}{n-2}-\frac{4}{5}\right\}$ \\ \hline
$\hwedge$ & $\ts\frac{1}{12}\frac{(v)_3}{(n)_3}$ & $\ts +\frac{(n-v)v}{(n)_2}$      & $\ts -\frac{2(n-v)(v)_2}{3(n)_3}$ & $\ts +\frac{5(n-v)v}{3(n)_2}\max\left\{0, \frac{v-1}{n-2}-\frac{9}{10}\right\}$ \\ \hline
$\hpar$ & $\ts\frac{1}{3}\frac{(v)_3}{(n)_3}$    & $\ts +\frac{(n-v)v}{(n)_2}$      & $\ts -\frac{2(n-v)(v)_2}{3(n)_3}$ & \\ \hline
$\hw$ & $\ts -\frac{3}{4}\frac{(v)_3}{(n)_3}$    & $\ts +\frac{(n-v)v}{(n)_2}$      & $\ts -\frac{(n-v)(v)_2}{2(n)_3}$ & \\ \hline
$\hpath$ & & $\ts \frac{(n-v)v}{(n)_2}$          & $\ts -\frac{(n-v)(v)_2}{2(n)_3}$ & \\ \hline
$\htri$ & $\ts\frac{1}{4}\frac{(v)_3}{(n)_3}$    & $\ts +\frac{(n-v)v}{(n)_2}$      & $\ts -\frac{(n-v)(v)_2}{2(n)_3}$ & \\ \hline
$\hcyc$ & $\ts -\frac{(v)_3}{(n)_3}$             & $\ts +\frac{(n-v)v}{(n)_2}$      & $\ts -\frac{(n-v)(v)_2}{3(n)_3}$ & \\ \hline
$\hq$ & $\ts -\frac{1}{4}\frac{(v)_3}{(n)_3}$    & $\ts +\frac{(n-v)v}{(n)_2}$      & $\ts -\frac{(n-v)(v)_2}{3(n)_3}$ & \\ \hline
$\hfive$ & $\ts -\frac{1}{2}\frac{(v)_3}{(n)_3}$ & $\ts +\frac{(n-v)v}{(n)_2}$      & $\ts -\frac{(n-v)(v)_2}{6(n)_3}$ & \\ \hline
$\hsix$ & & $\ts \frac{(n-v)v}{(n)_2}$           & & \\ \hline
\end{tabular}
\caption{The values of $d^*(H)$ for the eleven distinct nonisomorphic graphs on 4 vertices.}
\label{table:dstar}
\end{table}

\begin{table}[ht]\centering
\begin{tabular}{||r|ll||}\hline
$H$      & \multicolumn{2}{l||}{$d^*(H)$} \\ \hline\hline
$\hzero$ & $\ts\frac{(n-v)v}{(n)_2}$  & $\ts +\frac{1}{4}\frac{(v)_4}{(n)_4} -\frac{(n-v)(v)_2}{(n)_3}$ \\ \hline
$\hone$ & $\ts\frac{(n-v)v}{(n)_2}$   & $\ts +\frac{1}{6}\frac{(v)_3}{(n)_3} -\frac{5(n-v)(v)_2}{6(n)_3}$ \\ \hline
$\hwedge$ & $\ts\frac{(n-v)v}{(n)_2}$ & $\ts +\frac{1}{12}\frac{(v)_3}{(n)_3} -\frac{2(n-v)(v)_2}{3(n)_3}$ \\ \hline
$\hpar$ & $\ts\frac{(n-v)v}{(n)_2}$   & $\ts +\frac{1}{3}\frac{(v)_3}{(n)_3} -\frac{2(n-v)(v)_2}{3(n)_3}$ \\ \hline
$\hw$ & $\ts\frac{(n-v)v}{(n)_2}$     & $\ts -\frac{3}{4}\frac{(v)_3}{(n)_3} -\frac{(n-v)(v)_2}{2(n)_3}$ \\ \hline
$\hpath$ & $\ts\frac{(n-v)v}{(n)_2}$  & $\ts -\frac{(n-v)(v)_2}{2(n)_3}$ \\ \hline
$\htri$ & $\ts\frac{(n-v)v}{(n)_2}$   & $\ts +\frac{1}{4}\frac{(v)_3}{(n)_3} -\frac{(n-v)(v)_2}{2(n)_3}$ \\ \hline
$\hcyc$ & $\ts\frac{(n-v)v}{(n)_2}$   & $\ts -\frac{(v)_3}{(n)_3} -\frac{(n-v)(v)_2}{3(n)_3}$ \\ \hline
$\hq$ & $\ts\frac{(n-v)v}{(n)_2}$     & $\ts -\frac{1}{4}\frac{(v)_3}{(n)_3} -\frac{(n-v)(v)_2}{3(n)_3}$ \\ \hline
$\hfive$ & $\ts\frac{(n-v)v}{(n)_2}$  & $\ts -\frac{1}{2}\frac{(v)_3}{(n)_3} -\frac{(n-v)(v)_2}{6(n)_3}$ \\ \hline
$\hsix$ & $\ts \frac{(n-v)v}{(n)_2}$  & \\ \hline
\end{tabular}
\caption{The values of $d^*(H)$ in the case of $4\leq v\leq (2n-1)/3$.}
\label{table:dplusc_one}
\end{table}

\begin{table}[ht]\centering
\begin{tabular}{||r|rcrcr||}\hline
$H$       & \multicolumn{5}{l||}{$\gamma c(H)$} \\ \hline\hline
$\hzero$  & $0$         & $+$ & $0$         & $=$ & $0$ \\ \hline
$\hone$   & $0$         & $+$ & $4\gamma$   & $=$ & $4\gamma$ \\ \hline
$\hwedge$ & $4\gamma$   & $+$ & $0$         & $=$ & $4\gamma$ \\ \hline
$\hpar$   & $-16\gamma$ & $+$ & $8\gamma$   & $=$ & $-8\gamma$ \\ \hline
$\hw$     & $0$         & $+$ & $0$         & $=$ & $0$ \\ \hline
$\hpath$  & $-4\gamma$  & $+$ & $0$         & $=$ & $-4\gamma$ \\ \hline
$\htri$   & $12\gamma$  & $+$ & $-12\gamma$ & $=$ & $0$ \\ \hline
$\hcyc$   & $0$         & $+$ & $0$         & $=$ & $0$ \\ \hline
$\hq$     & $0$         & $+$ & $-4\gamma$  & $=$ & $-4\gamma$ \\ \hline
$\hfive$  & $0$         & $+$ & $4\gamma$   & $=$ & $4\gamma$ \\ \hline
$\hsix$   & $0$         & $+$ & $24\gamma$  & $=$ & $24\gamma$ \\ \hline
\end{tabular}
\caption{The values of $\gamma c(H)$ for the eleven distinct nonisomorphic graphs on 4 vertices. The first term is the contribution from the expression (\ref{eq:sq:1}), the second term from the expression in (\ref{eq:sq:2}). Their sum is the last term.}
\label{table:see}
\end{table}

\begin{table}[ht]\centering
\begin{tabular}{||r|lllll||}\hline
$H$       & \multicolumn{5}{l||}{$d^*(H)+\gamma c(H)$} \\ \hline\hline
$\hzero$  & $\ts\frac{1}{4}\frac{(v)_4}{(n)_4}$   & $\ts +\frac{(n-v)v}{(n)_2}$ & $\ts -\frac{(n-v)(v)_2}{(n)_3}$   & $\ts +\frac{3(n-v)v}{(n)_2}\max\left\{0, \frac{v-1}{n-2}-\frac{2}{3}\right\}$ & \\ \hline
$\hone$   & $\ts\frac{1}{6}\frac{(v)_3}{(n)_3}$   & $\ts +\frac{(n-v)v}{(n)_2}$ & $\ts -\frac{5(n-v)(v)_2}{6(n)_3}$ & $\ts +\frac{5(n-v)v}{2(n)_2}\max\left\{0, \frac{v-1}{n-2}-\frac{4}{5}\right\}$ & $\ts +4\left(\frac{1}{96}-\frac{(n-v)v}{24(n)_2}\right)$ \\ \hline
$\hwedge$ & $\ts\frac{1}{12}\frac{(v)_3}{(n)_3}$  & $\ts +\frac{(n-v)v}{(n)_2}$ & $\ts -\frac{2(n-v)(v)_2}{3(n)_3}$ & $\ts +\frac{5(n-v)v}{3(n)_2}\max\left\{0, \frac{v-1}{n-2}-\frac{9}{10}\right\}$ & $\ts +4\left(\frac{1}{96}-\frac{(n-v)v}{24(n)_2}\right)$ \\ \hline
$\hpar$   & $\ts\frac{1}{3}\frac{(v)_3}{(n)_3}$   & $\ts +\frac{(n-v)v}{(n)_2}$ & $\ts -\frac{2(n-v)(v)_2}{3(n)_3}$ & & $\ts -8\left(\frac{1}{96}-\frac{(n-v)v}{24(n)_2}\right)$ \\ \hline
$\hw$     & $\ts -\frac{3}{4}\frac{(v)_3}{(n)_3}$ & $\ts +\frac{(n-v)v}{(n)_2}$ & $\ts -\frac{(n-v)(v)_2}{2(n)_3}$  & & \\ \hline
$\hpath$  &                                       & $\ts \frac{(n-v)v}{(n)_2}$  & $\ts -\frac{(n-v)(v)_2}{2(n)_3}$  & & $\ts -4\left(\frac{1}{96}-\frac{(n-v)v}{24(n)_2}\right)$ \\ \hline
$\htri$   & $\ts\frac{1}{4}\frac{(v)_3}{(n)_3}$   & $\ts +\frac{(n-v)v}{(n)_2}$ & $\ts -\frac{(n-v)(v)_2}{2(n)_3}$  & & \\ \hline
$\hcyc$   & $\ts -\frac{(v)_3}{(n)_3}$ & $\ts +\frac{(n-v)v}{(n)_2}$ & $\ts -\frac{(n-v)(v)_2}{3(n)_3}$ & & \\ \hline
$\hq$     & $\ts -\frac{1}{4}\frac{(v)_3}{(n)_3}$ & $\ts +\frac{(n-v)v}{(n)_2}$ & $\ts -\frac{(n-v)(v)_2}{3(n)_3}$  & & $\ts -4\left(\frac{1}{96}-\frac{(n-v)v}{24(n)_2}\right)$ \\ \hline
$\hfive$  & $\ts -\frac{1}{2}\frac{(v)_3}{(n)_3}$ & $\ts +\frac{(n-v)v}{(n)_2}$ & $\ts -\frac{(n-v)(v)_2}{6(n)_3}$  & & $\ts +4\left(\frac{1}{96}-\frac{(n-v)v}{24(n)_2}\right)$ \\ \hline
$\hsix$   &                                       & $\ts \frac{(n-v)v}{(n)_2}$  &                                   & & $\ts +24\left(\frac{1}{96}-\frac{(n-v)v}{24(n)_2}\right)$ \\ \hline
\end{tabular}
\caption{The values of $d^*(H)+\gamma c(H)$ in the case of $v\geq 2n/3$.}
\label{table:dplusc_twoa}
\end{table}

\begin{table}[ht]\centering
\begin{tabular}{||r|lll||}\hline
$H$       & \multicolumn{3}{l||}{$d^*(H)+\gamma c(H)-1/4$} \\ \hline\hline
$\hzero$  & $\ts -\frac{1}{4}$  & $\ts +\frac{1}{4}\frac{(v)_4}{(n)_4}$  & $\ts +\frac{(n-v)v(2v-n)}{(n)_3}$ \\ \hline
$\hone$   & $\ts -\frac{5}{24}$ & $\ts +\frac{1}{6}\frac{(v)_3}{(n)_3}$  & $\ts +\frac{(n-v)v}{6(n)_3}\max\{5n-5v-5, 10v-7n+4\}$ \\ \hline
$\hwedge$ & $\ts -\frac{5}{24}$ & $\ts +\frac{1}{12}\frac{(v)_3}{(n)_3}$ & $\ts +\frac{(n-v)v}{6(n)_3}\max\{5n-4v-6, 6v-4n+2\}$ \\ \hline
$\hpar$   & $\ts -\frac{1}{3}$  & $\ts+\frac{1}{3}\frac{(v)_3}{(n)_3}$   & $\ts +\frac{(n-v)v(4n-2v-6)}{3(n)_3}$ \\ \hline
$\hw$     & $\ts -\frac{1}{4}$  & $\ts -\frac{3}{4}\frac{(v)_3}{(n)_3}$  & $\ts +\frac{(n-v)v(2n-v-3)}{2(n)_3}$ \\ \hline
$\hpath$  & $\ts -\frac{7}{24}$ &                                        & $\ts +\frac{(n-v)v(7n-3v-11)}{6(n)_3}$ \\ \hline
$\htri$   & $\ts -\frac{1}{4}$  & $\ts +\frac{1}{4}\frac{(v)_3}{(n)_3}$  & $\ts +\frac{(n-v)v(2n-v-3)}{2(n)_2}$ \\ \hline
$\hcyc$   & $\ts -\frac{1}{4}$  & $\ts -\frac{(v)_3}{(n)_3}$ & $\ts +\frac{(n-v)v(3n-v-5)}{3(n)_3}$ \\ \hline
$\hq$     & $\ts -\frac{7}{24}$ & $\ts -\frac{1}{4}\frac{(v)_3}{(n)_3}$  & $\ts +\frac{(n-v)v(7n-2v-12)}{6(n)_3}$ \\ \hline
$\hfive$  & $\ts -\frac{5}{24}$ & $\ts -\frac{1}{2}\frac{(v)_3}{(n)_3}$  & $\ts +\frac{(n-v)v(5n-v-9)}{6(n)_3}$ \\ \hline
$\hsix$   & $0$ & & \\ \hline
\end{tabular}
\caption{Simplified values of $d^*(H)+\gamma c(H)-1/4$ in the case of $v\geq 2n/3$.}
\label{table:dplusc_twob}
\end{table}

\begin{table}[ht]\centering
\begin{tabular}{||r|l|rcr||}\hline
$H$       & $g_H(x)$ & \multicolumn{3}{r||}{max over $\frac{2}{3}\leq x\leq 1$} \\ \hline\hline
$\hone$   & $\ts -\frac{5}{24}+\frac{x^3}{6}+\frac{(1-x)x}{6}\max\{5-5x,10x-7\}$ & $\ts -\frac{1}{24}$    & @ & $x=1$ \\ \hline
$\hwedge$ & $\ts -\frac{5}{24}+\frac{x^3}{12}+\frac{(1-x)x}{6}\max\{5-4x,6x-4\}$ & $\ts -\frac{7}{72}$    & @ & $\ts x=\frac{2}{3}$ \\ \hline
$\hpar$   & $\ts -\frac{1}{3}+\frac{x^3}{3}+\frac{(1-x)x(4-2x)}{3}$              & $\ts 0$                & @ & $x=1$ \\ \hline
$\hw$     & $\ts -\frac{1}{4}-\frac{3x^3}{4}+\frac{(1-x)x(2-x)}{2}$              & $\ts -\frac{35}{108}$  & @ & $\ts x=\frac{2}{3}$ \\ \hline
$\hpath$  & $\ts -\frac{7}{24}+\frac{(1-x)x(7-3x)}{6}$                           & $\ts -\frac{23}{216}$  & @ & $\ts x=\frac{2}{3}$ \\ \hline
$\htri$   & $\ts -\frac{1}{4}+\frac{x^3}{4}+\frac{(1-x)x(2-x)}{2}$               & $\ts 0$                & @ & $\ts x=1$ \\ \hline
$\hcyc$   & $\ts -\frac{1}{4}-x^3+\frac{(1-x)x(3-x)}{3}$                         & $\ts -\frac{121}{324}$ & @ & $\ts x=\frac{2}{3}$ \\ \hline
$\hq$     & $\ts -\frac{7}{24}-\frac{x^3}{4}+\frac{(1-x)x(7-2x)}{6}$             & $\ts -\frac{101}{648}$ & @ & $\ts x=\frac{2}{3}$ \\ \hline
$\hfive$  & $\ts -\frac{7}{24}-\frac{x^3}{2}+\frac{(1-x)x(5-x)}{6}$              & $\ts -\frac{127}{648}$ & @ & $\ts x=\frac{2}{3}$ \\ \hline
\end{tabular}
\caption{The $g_H(x)$ functions that determine the asymptotic value of $d^*(H)+\gamma c(H)-1/4$.}
\label{table:xfunct}
\end{table}

\FloatBarrier
\end{document}

%% file: diamond-macro.tex
\usepackage{amsfonts,enumerate,placeins,amsmath}
\newtheorem{thm}{Theorem}
\newtheorem{conj}[thm]{Conjecture}
\newtheorem{lem}[thm]{Lemma}
\newtheorem{defn}[thm]{Definition}
\newtheorem{cl}[thm]{Claim}

\newtheorem{notat}[thm]{Notation}

\newcommand{\ts}{\textstyle}
\DeclareMathOperator{\La}{La}
\DeclareMathOperator{\Lam}{\Lambda}
\DeclareMathOperator{\Lamstar}{\Lambda^*}
\DeclareMathOperator{\rmdef}{def}
\DeclareMathOperator{\rmmin}{min}
\newcommand{\PP}{\mathcal{P}}
\newcommand{\F}{\mathcal{F}}

\newcommand{\T}{\mathcal{T}}
\newcommand{\U}{\mathcal{U}}
\newcommand{\nntwo}{{n \choose \left\lfloor n/2 \right\rfloor}}
\newcommand{\booltwo}{\mathcal{B}_2}
\newcommand{\bool}{\mathcal{B}}

\newenvironment{pfct}[1]{\noindent{\bf Proof of #1.\,}}{\hfill$\Box$ \\}

\newcommand{\booln}{\mathcal{B}_n}

\newcommand{\hzero}{H_0}
\newcommand{\betazero}{\beta_0}

\newcommand{\hone}{H_1}

\newcommand{\hwedge}{H_{\wedge}}

\newcommand{\hpar}{H_{\|}}

\newcommand{\htri}{H_{\triangle}}

\newcommand{\hw}{H_{\dashv}}

\newcommand{\hpath}{H_{\sqcup}}

\newcommand{\hcyc}{H_{\square}}

\newcommand{\hq}{H_{Q}}

\newcommand{\hfive}{H_5}

\newcommand{\hsix}{H_6}

\newcommand{\ove}{\overline{e}}